\documentclass{tran-l}
\begin{document}

\newcommand{\dt}{\displaystyle}

\title[THE ALGORITHM FOR SOLVING OF THE GRAPH ISOMORPHISM
PROBLEM]
 {THE DIRECT ALGORITHM FOR SOLVING OF THE GRAPH ISOMORPHISM
PROBLEM}

\author{ Rashit T.~Faizullin, Alexander
V.~Prolubnikov}

\address{Omsk
State University, Omsk, Russia}

\email{faizulin@univer.omsk.su,prolubnikov@univer.omsk.su}

\maketitle

\begin{abstract}
{\small We propose an algorithm for solving of the graph
isomorphism problem. Also, we introduce the new class of graphs
for which the graph isomorphism problem can be solved polynomially
using the algorithm.}
\end{abstract}

\section*{\bf INTRODUCTION}

Complexity of the graph isomorphism problem is not definite still.
The graph iso\-mor\-phism problem can not be placed in any of
known complexity classes as it was stated in \cite{1}. It is not
proved that the problem is $NP$-complete and there is no
polynomial-time algorithm for solving of the problem too.
Polynomial-time algorithms exists only for the special cases of
the graph isomorphism problem \cite{2, 3, 4}.

We propose an algorithm which is based on the methods of linear
algebra. The algorithm gives a solution of the problem for graphs
from the class that we define in this paper. Definition of the
class is closely related with graph isomorphism notion. We checks
weather graphs belongs to the class or not at the process of the
algorithm operating. We couldn't find a counter example for the
algorithm, i.e., we couldn't find the pair of isomorphic graphs
which can't be defined as isomorphic by the algorithm. In
particular, at the process of numerous numerical experiments, we
couldn't find a counter example too. Also, such graphs as regular
graphs, which traditionally gives hardest cases of the graph
isomorphism problem, they belongs to this class
as numerical experiments shows. 

At the process of the algorithm operating, we implement
perturbations of modified adjacency matrices of the graphs and
solve the systems of linear equations associated with them. The
algorithm is direct, i.e., the algorithm is not a variation of
backtracking scheme. There is no search tree, which growth may
become uncontrolled at the process of the algorithm operation. The
isomorphism of the graphs checking at most at $n$ iterations of
the algorithm, where $n$ is a number of vertices of graphs. If the
graphs belong to the class and they are isomorphic, then one
isomorphism is presented as a result. It is shown that complexity
of the presented algorithm is equal to $O(n^6\log n)$ in sense of
using elementary operations.

The paper is organized as follows. In $\S 1$, we presents the
basic scheme of the algorithm and gives a theoretical basis of it
using terms of linear algebra. In $\S 2$, we propose the numerical
example of the algorithm implementing. In $\S 3$, we consider the
main computational problems the algorithm deal with; its overall
complexity is considered.

\section{\bf THE BASE SCHEME OF THE ALGORITHM}

\noindent {\bf The graph isomorphism problem. Formulation 1.}
Suppose $G_A=\langle V_A,E_A\rangle $ and $G_B=\langle
V_B,E_B\rangle $ are two non-weighted non-oriented graphs, where
$V_A,V_B$ are sets of their vertices and $E_A,E_B$ are sets of
their edges. $|V_A|=|V_B|,\ |E_A|=|E_B|$.

Graphs $G_A=\langle V_A,E_A\rangle$ and $G_B=\langle
V_B,E_B\rangle $ are said to be isomorphic if there exists a
bijection $\varphi : V_B\to V_A$ such that $(i,j)\in E_B
\Leftrightarrow (\varphi(i),\varphi(j))\in E_A$. One must find
this $\varphi$ or prove that there is no such bijection. \qed

Let us denote isomorphism from graph $G_A$ to $G_B$ as $G_A\simeq
G_B$. Let $A_0$ be an adjacency matrix of graph $G_A$, i.e.,
$A_0=(a_{ij}^0)$, where
$$a_{ij}^0= \left\{
\begin{array}{ll}
1,&\mbox{ if } (i,j)\in E_A,\\
0,&\mbox{ else}.\\
\end{array}
\right.
$$
Let $B_0$ be an adjacency matrix of $G_B$.

We can uniquely set a permutation matrix $P_\varphi=(p_{ij})$ in
correspondence for every bijection $\varphi : V_B\to V_A$.
Elements of corresponded permutation matrix $P$ are followed:
$$
p_{ij}=\left\{
\begin{array}{ll}
1,&\mbox{ if } i=\varphi(j),\\
0,&\mbox{ else}.\\
\end{array}
\right. $$

The following formulation of the graph isomorphism problem is
equivalent to the first one.

\noindent {\bf The graph isomorphism problem. Formulation 2.}
Suppose $G_A=\langle V_A,E_A\rangle $ and $G_B=\langle
V_B,E_B\rangle $ are two non-weighted non-oriented graphs.
$|V_A|=|V_B|,\ |E_A|=$ $=|E_B|$. Let $A_0$ and $B_0$ be adjacency
matrices of the graphs $G_A$ and $G_B$.

One must find a permutation matrix $P$ such that $A_0=PB_0P^{-1}$
or prove that there is no such permutation matrix. \qed

Let $A_0$ be an adjacency matrix of $G_A=\langle V_A,E_A\rangle $.
Let $D_{A_0}$ be a diagonal matrix of the following form: $$\left
(
\begin{array}{cccc}
  d_{11}^A  &   0   &\ldots&   0\\
   0   &   d_{22}^A &\ldots&   0\\
 \vdots& \vdots&\ddots&\vdots\\
   0   &   0   &\ldots&   d_{nn}^A
\end{array}
\right),$$ where $$d_{ii}^A=\sum\limits_{j=1}^na_{ij}^0+d=d_i+d$$
and $d$ is a maximal degree of $G_A$ vertices, $d_i$ is a degree
of vertex $i\in V_A$. Matrix $D_{B_0}$ is constructed in similar
manner according to matrix $B_0$.

The algorithm operates with graph matrices of the following form
$$A=A_0+D_{A_0},\ B=B_0+D_{B_0}. \eqno (1)$$ $A$ and $B$ are positive definite matrices with
diagonal predominance. Let us call graph matrices of the form (1)
as graph matrices.

It is clear that the following formulation of the graph
isomorphism problem is equivalent to the second one and hence it
is equivalent to the first one.

\noindent {\bf The graph isomorphism problem. Formulation 3.}
Suppose $G_A=\langle V_A,E_A\rangle $ and $G_B=\langle
V_B,E_B\rangle $ are two non-weighted non-oriented graphs.
$|V_A|=|V_B|,\ |E_A|=$ $=|E_B|$. Let $A$ and $B$ be graph matrices
of the form (1) for $G_A$ and $G_B$.

One must find a permutation matrix $P$ such that $A=PBP^{-1}$ or
prove that there is no such permutation matrix. \qed

Formulation 3 of the graph isomorphism problem may be considered
as a reduction of the initial problem to the problem of
isomorphism checking of the same graphs with additional weighted
loops. Matrices of the form (1) are adjacency matrices of these
graphs.

Consider the following systems of linear algebraic equations.
$$Ax=e_j,\ By=e_k, \eqno (2)$$
where $\{e_j\}_{j=1}^n$ is a standard basis in $\mathbb{R}^n$ and
$A$ and $B$ are matrices of the form (1) (graph matrices) for
graphs $G_A$ and $G_B$. Since $A$ and $B$ are positive definite
matrices, every system of equations in (2) has solution and the
solution is unique. Let $x_j$ be a solution of the system of
equations $Ax=e_j$ and let $y_k$ be solution of the system of
equations $By=e_k$. Vector $x_j=(x_{j1},\ldots, x_{jn})$ has the
following components
$$x_{ji}=\frac{A_{ij}}{\det A},$$ where $A_{ij}$ is a cofactor of
element $a_{ij}$ of matrix $A$, i.e. $x_j$ is $j$-th column of
matrix $A^{-1}$.

Note that $A=PBP^{-1}$ iff $A^{-1}=PB^{-1}P^{-1}$.

\bigskip

Put $A^k=A+\sum\limits_{j=1}^k\varepsilon_jE^j,$ where $$E^j=
\left (
\begin{array}{ccccc}
0      &\ldots &   0   &\ldots & 0\\ \vdots &\ddots &\vdots & &
\vdots\\ 0      &\ldots &   1   &\ldots & 0\\ \vdots &
&\vdots &\ddots & \vdots\\ 0      &\ldots &   0   &\ldots & 0\\
\end{array}
\right )
\begin{array}{cc}
- & 1\\ \vdots\\ - & j\\ \vdots\\ - & n
\end{array}
$$ i.e., $E^j$ is a diagonal matrix with the only nonzero element. It is an
element at its $j$-th row and it is equal to 1. $\varepsilon_j\in
\mathbb{R},\ \varepsilon_j>0$.

Let $G_A^{1,\ldots ,k}$ be the graph which adjacency matrix is
matrix $A^k$.

Put $B^l=B+\sum\limits_{j=1}^l\varepsilon_jE^{k_j}$ and let
$G_B^{k_1,\ldots ,k_l}$ be the graph which adjacency matrix is
matrix $B^l$. Note that $B^l$ is defined by sequence $\{k_1,\ldots
,k_l\}$ because we perturb the diagonal elements of $B$ by
$\varepsilon_j$ accordingly this sequence while $A^k$ is always
defined by the only sequence $\{1,\ldots ,k\}$.

{\bf Definition 1.} We say that graphs $G_A=\langle
V_A,E_A\rangle$ and $G_B=\langle V_A,E_A\rangle$ are {\it similar}
and denote it as $G_A\sim G_B$ whenever there exists bijections
$\varphi,\varphi_i:V_B\to V_A$ such that
$x_{ij}=y_{\varphi(i)\varphi_i(j)}$, $(i,j=\overline{1,n})$.\qed

{\bf Definition 2.} We say that graph $G_A$ is {\it
reconstructible} whenever   $$G_A\sim G_B \mbox{ iff } G_A\simeq
G_B$$ and it holds for every graph $G_B$.\qed

Let $G_A^{ij}$ ($i,j=\overline{1,n}$) be the graph which adjacency
matrix is $A_{(i,j)}$, where $A_{(i,j)}$ is obtained from matrix
$A$ of the form (1) by deleting its $i$-th row and $j$-th column.
If $i\neq j$, then matrix $A_{(i,j)}$ is not symmetrical at
general case. $G_A^{ij}$ is a weighted oriented graph with loops.
$G_A^{ii}$ ($i=\overline{1,n}$) is a graph which is obtained from
$G_A$ by deleting vertex $i\in V_A$ and edges which are incident
to $i$. Let us call the graph $G_A^{ij}$, ($i,j=\overline{1,n}$)
as {\it associated graph}.

Let $$\Lambda_A^{ij}=\prod\limits_{k=1}^{n-1}\lambda_k
(G_A^{ij}),\ \Lambda_B^{ij}=\prod\limits_{k=1}^{n-1}\lambda_k
(G_B^{ij}),$$
$$\Lambda_A=\prod\limits_{k=1}^n\lambda_k (G_A),\
\Lambda_B=\prod\limits_{k=1}^n\lambda_k (G_B),$$ where $\lambda_k
(G_A^{ij})$ is $k$-th eigenvalue of matrix $A_{(i,j)}$, $\lambda_k
(G_B^{ij})$ is $k$-th eigenvalue of matrix $B_{(i,j)}$, $\lambda_k
(G_A)$ is $k$-th eigenvalue of matrix $A$, $\lambda_k (G_B)$ is
$k$-th eigenvalue of matrix $B$. Thus
$$x_{ij}=\frac{A_{ij}}{\det A}=\frac{\Lambda_A^{ij}}{\Lambda_A},
y_{ij}=\frac{B_{ij}}{\det B}=\frac{\Lambda_B^{ij}}{\Lambda_B}.$$
This means that graphs $G_A$ and $G_B$ are similar iff there
exists bijections $\varphi,\varphi_i:V_B\to V_A$ such that for
graphs $G_A^{ij}$ and $G_B^{\varphi(i)\varphi_i(j)}$
$(i,j=\overline{1,n}$) the following holds:
$$\prod\limits_{k=1}^{n-1}\lambda_k
(G_A^{ij})=\prod\limits_{k=1}^{n-1}\lambda_k
(G_B^{\varphi(i)\varphi_i(j)}).$$

There are a lot of ways for constructing cospectral non-isomorphic
graphs. Different types of graph matrices may be considered. But
note that there is no universal way for this construction. It
depends on both the graph representing matrices and graphs
structure. Checking similarity of the graphs, we doesn't check
their cospectrality. Indeed, we check more slight spectral
property of the graphs, namely we check the equality of products
of eigenvalues from corresponded associated graphs spectrums. But
we checks coincidence of this characteristics for ample quantity
of associated graphs at one time. There are $n(n-1)/2$ of
associated graphs for one graph. Checking this characteristics at
every iteration of the algorithm gives us a powerful method for
solving of the graph isomorphism problem. We can not find a
counter example for the algorithm even though there are numerous
explorations have been implemented. In particular, there is no
non-reconstructible graphs in the sense of definition 2 at the
library which is located at the url {\it
http://amalfi.dis.unina.it/graph/}. This library is used for
testing of the most efficient algorithms which are designed for
solving of the problem \cite{5}.

If $G_A\simeq G_B$ and $\varphi$ is an isomorphism from $G_A$ to
$G_B$, then we have

$\left (
\begin{array}{ccccc}
  x_{11}  & \ldots & x_{1j} & \ldots &   x_{1n}\\
  \vdots  & \ddots & \vdots &        &   \vdots\\
  x_{j1}  & \ldots & x_{jj} & \ldots &   x_{jn}\\
  \vdots  &        & \vdots & \ddots &   \vdots\\
  x_{n1}  & \ldots & x_{nj} & \ldots &   x_{nn}
\end{array}
\right)=$\hfill $$=\left (
\begin{array}{ccccc}
  y_{\varphi (1)\varphi (1)} & \ldots & y_{\varphi (1)\varphi (j)} & \ldots &   y_{\varphi (1)\varphi (n)}\\
  \vdots          & \ddots & \vdots          &        &   \vdots\\
  y_{\varphi (j)\varphi (1)} & \ldots & y_{\varphi (j)\varphi (j)} & \ldots &   y_{\varphi (j)\varphi (n)}\\
  \vdots          &        & \vdots          & \ddots &   \vdots\\
  y_{\varphi (n)\varphi (1)} & \ldots & y_{\varphi (n)\varphi (j)} & \ldots &   y_{\varphi (n)\varphi (n)}
\end{array}
\right).$$

{\bf Lemma 1.} Let $A=PBP^{-1}$ holds for matrices $A$ and $B$,
where $P$ is a permutation matrix. If there are no columns of
matrix $A^{-1}$ such that they are coincide to within the
permutation of its components, then $P$ is the only permutation
matrix for which it holds.

{\bf Proof.} The proof of the lemma is obvious. \qed

{\bf Remark 1.} Clearly, if we take the matrix $A$ itself as a
matrix $B$ in the conditions of Lemma 1, then we obtain the
following statement. If there is no columns of matrix $A^{-1}$
such that they are coincide to within the permutation of its
components, then $\Gamma (G_A)$ is trivial (or $|\Gamma
(G_A)|=1$).

Presented algorithm operates with perturbed matrices which are
obtained from matrices of the form (1). The perturbations
implements by adding matrices of the form $\varepsilon E^j$ to the
graph matrices. Let $A^0=A,\ B^0=B.$ At one iteration of the
algorithm, we consider one matrix $A^{j}$ and matrices $B_k^{j}$,
where
$$A^{j}=A^{j-1}+\varepsilon E^j,\ B_k^{j}=B^{j-1}+\varepsilon
E^k.$$ Varying $k$, we finds the vertex $k_j\in V_B$ and we set
$k_j\in V_B$ in correspondence to $j\in V_A$ iff $G_A^{1,\ldots
,j}\sim G_B^{k_1,\ldots ,k_{j-1},k_j}$. Perturbations performs
until there are no columns of $A^{-1}$ that are coinciding within
the permutation of its components, i.e., perturbations performs
until automorphisms groups of the graphs become trivial. Let
"$j\leftrightarrow k_j$" denote the correspondence that we set at
the process of the algorithm implementing.

\begin{center}
{{\bf The algorithm of graph isomorphism testing}\\
Base scheme}
\end{center}

\noindent Step 0. $j:=1,\ k_1:=0$. Go to step 1.\\ \noindent Step
1. If $j>n$, then go to step 7 else chose $\varepsilon_k$, $l:=1$,
go to step 2.\\ \noindent Step 2. If $l>n$, then go to step 6 else
if $l=k_s$ ($1\le s\le j-1$), then go to step 4 else go to step
3.\\ \noindent Step 3. If $G_A^{1,\ldots ,j}\sim G_B^{k_1,\ldots
,k_{j-1},l}$,
then $k_j:=l$ and go to step 5 else go to step 4.\\
\noindent Step 4. $l:=l+1$. Go to step 2.\\ \noindent Step 5.
$j:=j+1$. Go to step 1.\\ \noindent Step 6. Graphs $G_A$ and $G_B$
are not isomorphic or they are not reconstructible. Stop the algorithm.\\
\noindent Step 7. Graphs $G_A$ and $G_B$ are reconstructible and
isomorphic. Let $\varphi: V_B\to V_A$ be a bijection such that
$\varphi (k_j)=j$. $\varphi$ is an isomorphism from $G_A$ to
$G_B$. Stop the algorithm.

{\bf Lemma 2.} Let $P$ be a permutation matrix and let $\varphi :
V_B \to V_A$ be a bijection that is corresponded to $P$.  Let
$A^j$ and $B^j$ be the matrices which are obtained from $A$ and
$B$ by perturbations of its diagonal elements accordingly to
sequences $\{j\}_{j=1}^n$ and $\{k_j\}_{j=1}^n$, where $\varphi
(k_j)=j$.

Then $$A=PBP^{-1}\Leftrightarrow \forall j\ A^j=PB^jP^{-1}.$$

{\bf Proof.} $A^j=A^{j-1}+\varepsilon_jE^j,\
B^j=B^{j-1}+\varepsilon_jE^{k_j}.$ Since $$\varphi
(k_j)=j\Leftrightarrow E^j=PE^{k_j}P^{-1},$$ it follows that
$$A^j=PB^jP^{-1}\Leftrightarrow
A^{j-1}+\varepsilon_jE^j=$$
$$=P(B^{j-1}+\varepsilon_jE^{k_j})P^{-1}=PB^{j-1}P^{-1}+\varepsilon_jPE^{k_j}P^{-1}\Leftrightarrow $$
$$\Leftrightarrow (A^{j-1}=PB^{j-1}P^{-1}\mbox{ and }j=\varphi (k_j)).$$
The last is true whenever it holds for $j=1$, i.e.,
$A=PBP^{-1}$.\qed

It follows from lemma 2 that $G_A\simeq G_B$ iff there exists a
permutation matrix $P$ such that $A^j=PB^jP^{-1}$ holds for
$j=\overline{1,n}$ whenever $k_j=\varphi(j)$. I.e., if $G_A\simeq
G_B$, then, at every $j$-th iteration, the following must be true:
if $A^j=A^{j-1}+\varepsilon_jE^j$, then
$A^j=PB^jP^{-1}\Leftrightarrow
B^j=B^{j-1}+\varepsilon_jE^{\varphi^{-1} (j)}.$

By definition, put
$$R(A)=\{x_j\}_{j=1}^n,\
R(B)=\{y_k\}_{k=1}^n.$$

{\bf Definition 3.} A set $R(A)$ ($R(B)$) is called {\it simple}
if $\forall j\ \forall P\ \not\exists l\ (l\neq j): x_j=Px_l$
$(\forall k\ \forall P\not\exists i\ (k\neq i) \ y_k=Py_i),$ where
$P$ is a permutation matrix.\qed

Note that since graphs matrices have diagonal predominance, the
equality $x_j=Px_l$ implies that $x_{jj}=x_{ll}$ whenever $P$
corresponds to some isomorphism $\varphi$.

At the process of operating of the algorithm, it is necessary to
obtain matrices $A^n$ and $B^n$ such that $A^n=PB^nP^{-1}$ and
sets $R(A^n)$ and $R(B^n)$ are simple. Thus $P$ will be the only
permutation matrix that sets the isomorphisms of the graphs that
are corresponds to the matrices $A^n$ and $B^n$. For the graphs
that are corresponds to the matrices $A$ and $B$, this matrix will
be such permutation matrix too. Recall that, for checking the
isomorphism of reconstructible graphs, we may check the
coincidence of vectors $x_j$ and $y_k$ for correspondence within
to the permutation of their components.

Elements of inverse matrix $A^{-1}$ are continuous functions of
elements of matrix $A$. If $A'=A+\varepsilon_j E^j$ and
$B'=B+\varepsilon_j E^k$, then
$$\det A'=\det A+\varepsilon_jA_{jj},\
\det B'=\det B+\varepsilon_jB_{kk},\eqno (3)$$ Note that $\det
 A',\det A,A_{jj},\det B',\det B ,B_{kk}$ are positive. Again, let $P$ be a permutation matrix and let $\varphi$
be a bijection that is corresponded to $P$.

{\bf Lemma 3.} Let $$A=A^{j-1},\ B=B^{j-1},$$
$$A'=A+\varepsilon E^j,\ B_k'=B+\varepsilon E^k\
(k=\overline{1,n})$$ and let $x_j$ be solution of the system of
linear equations $A'x=e_j$ and let $y_k$ be solution of the system
of linear equations $B_k'y=e_k$.

If $A=PBP^{-1}$ holds for some permutation matrix $P$, then
$\exists !k: x_j=Py_k$.

{\bf Proof.} Let $x_i^0$ and $y_l^0$ be solutions of the followed
systems of linear equations:
$$Ax=e_i,\ By=e_l.$$
Then
$$x_i^0=\frac{1}{\det A}(A_{1i},\ldots,A_{ni}),\
y_l^0=\frac{1}{\det B}(B_{1l},\ldots,A_{nl}).$$ Since element
$a_{jj}$ of $A$ and element $b_{kk}$ of $B$ are the only elements
that we change by perturbations, then $A_{jt}'=A_{jt}$,
${(B_k')}_{kt}=B_{kt}$ ($t=\overline{1,n}$, where ${(B_k')}_{kt}$
is a cofactor of element $b_{tk}$ of matrix $B_k'$. $A$ and $B$
are matrices with diagonal predominance. So if $j=\varphi(k)$,
then $x_{jj}^0=y_{kk}^0$ since $x_j^0=Py_k^0$ and $\varphi$ is an
isomorphism from $G_A$ to $G_B$, i.e., $A_{jj}=B_{kk}$. Since
$A=PBP^{-1}$, we have $\det A=\det B$. Taking into account (3), we
get $\det A'=\det B_k'$ for every $\varepsilon$. Hence if
$x_j^0=Py_k^0$, then $x_j=Py_k$ too because
$$x_j=\frac{\det A}{\det A'}\times
\frac{1}{\det A}(A_{1j},\ldots,A_{nj})=\frac{\det A}{\det
A'}x_j^0,$$
$$y_k=\frac{\det B}{\det B_k'}\times
\frac{1}{\det B}(B_{1k},\ldots,A_{nk})=\frac{\det B}{\det
B_k'}y_k^0.\eqno (4)$$

Let us prove that if columns $x_j$ and $y_k$ were not coinciding
to within the permutation $P$ of its components before the
implemented perturbations, then they won't be coinciding after the
perturbations too.

While permutation of matrix rows connected with the same
permutation of its columns occurs, matrix diagonal elements
transforms into transformed matrix diagonal elements. Therefore
$x_j$ and $y_k$ may coincide to within the permutation of its
components only if $x_{jj}=y_{kk}$. There are two possibilities of
it. First, $A_{jj}=B_{kk}$ before the implemented perturbations,
second, $A_{jj}\neq B_{kk}$ before the perturbation. At first
case, taking into account (3), we get $\det A'=\det B_k'$. Thus if
$x_j^0$ and $y_k^0$ are not coinciding to within the permutation
of its components, then $x_j$ and $y_k$ are not coinciding to
within the permutation of its components too since (4) holds. If
$A_{jj}\neq B_{kk}$, then
$$x_{jj}=\frac{A_{jj}}{\det A'}=\frac{A_{jj}}{\det A+\varepsilon
A_{jj}}\neq \frac{B_{kk}}{\det B+\varepsilon
B_{kk}}=\frac{B_{kk}}{\det B_k'}=y_{kk},$$ where $\varepsilon$ is
an arbitrary. Therefore $x_j$ and $y_k$ are not coinciding to
within the permutation of its components.

\qed

It follows from lemma 3 that the perturbations may be implemented
in the way such that, in iverse matrix of every graph matrix,
there will be no columns that are coincides to within the
permutation of its components after the perturbations. Therefore,
it follows from lemma 1 that, at the process of the algorithm
operating, we transforms the initial graphs into the graphs that
have trivial automorphism group. Similarity of the graphs is
preserving while perturbations implements. In case of
reconstructible graphs, it follows from lemma 2 that
correspondence $j\leftrightarrow k_j$ is an isomorphism, where
$j\leftrightarrow k_j$ is a correspondence that is settled at the
iterations of the algorithm.

\section{\bf EXAMPLE OF THE ALGORITHM IMPLEMENTING}

Let us consider example of the algorithm implementing. Let $G_A$
and $G_B$ be isomorphic tested graphs. For these graphs, matrices
of the form (1) are followed matrices:
$$A=\left (
\begin{array}{cccccc}
8     & 1     & 1     & 1     & 1     & 0\\
1     & 7     & 1     & 0     & 0     & 1\\
1     & 1     & 7     & 0     & 0     & 1\\
1     & 0     & 0     & 7     & 1     & 1\\
1     & 0     & 0     & 1     & 7     & 1\\
0     & 1     & 1     & 1     & 1     & 8\\
\end{array}
\right ),\ B=\left (
\begin{array}{cccccc}
8     & 0     & 1     & 1     & 1     & 1\\
0     & 8     & 1     & 1     & 1     & 1\\
1     & 1     & 7     & 1     & 0     & 0\\
1     & 1     & 1     & 7     & 0     & 0\\
1     & 1     & 0     & 0     & 7     & 1\\
1     & 1     & 0     & 0     & 1     & 7\\
\end{array}
\right ).$$

\newpage
Let machine numbers mantissas length is equals to the number of
graph vertices $(n=6)$. Then \\ $ (A^0)^{-1}=\left (
\begin{array}{cccccc}
8     & 1     & 1     & 1     & 1     & 0\\
1     & 7     & 1     & 0     & 0     & 1\\
1     & 1     & 7     & 0     & 0     & 1\\
1     & 0     & 0     & 7     & 1     & 1\\
1     & 0     & 0     & 1     & 7     & 1\\
0     & 1     & 1     & 1     & 1     & 8\\
\end{array}
\right )^{-1}=$ $$= \left (
\begin{array}{cccccc}
0.134    &-0.018    &-0.018    &-0.018    &-0.018    & 0.008929\\
-0.018   & 0.15     &-0.016    &0.004464  &0.004464  &-0.018 \\
-0.018   &-0.016    & 0.15     &0.004464  &0.004464  &-0.018\\
-0.018   &0.004464  &0.004464  & 0.15     &-0.016    &-0.018\\
-0.018   &0.004464  &0.004464  &-0.016    & 0.15     &-0.018\\
0.008929 &-0.018    &-0.018    &-0.018    &-0.018    & 0.134\\
\end{array}
\right ),$$ and \\ $ (B^0)^{-1}=\left (
\begin{array}{cccccc}
8     & 0     & 1     & 1     & 1     & 1\\
0     & 8     & 1     & 1     & 1     & 1\\
1     & 1     & 7     & 1     & 0     & 0\\
1     & 1     & 1     & 7     & 0     & 0\\
1     & 1     & 0     & 0     & 7     & 1\\
1     & 1     & 0     & 0     & 1     & 7\\
\end{array}
\right )^{-1}=$
$$=
\left (
\begin{array}{cccccc}
0.134    &0.008929  &-0.018    &-0.018    &-0.018    &-0.018\\
0.008929 &0.134     &-0.018    &-0.018    &-0.018    &-0.018\\
-0.018   &-0.018    & 0.15     &-0.016    &0.004464  &0.004464\\
-0.018   &-0.018    &-0.016    & 0.15     &0.004464  &0.004464\\
-0.018   &-0.018    &0.004464  &0.004464  & 0.15     &-0.016\\
-0.018   &-0.018    &0.004464  &0.004464  &-0.016    & 0.15\\
\end{array}
\right ).$$


In this case, possible solutions of the graph isomorphism problem
are permutation matrices which are corresponded to the following
substitutions: $$ \left (
\begin{array}{cccccc}
1     & 2     & 3     & 4     & 5     & 6\\
1     & 3     & 4     & 5     & 6     & 2\\
\end{array}
\right ),\ \left (
\begin{array}{cccccc}
1     & 2     & 3     & 4     & 5     & 6\\
1     & 3     & 4     & 6     & 5     & 2\\
\end{array}
\right ),\ \left (
\begin{array}{cccccc}
1     & 2     & 3     & 4     & 5     & 6\\
1     & 4     & 3     & 5     & 6     & 2\\
\end{array}
\right ),\
$$
$$
\left (
\begin{array}{cccccc}
1     & 2     & 3     & 4     & 5     & 6\\
1     & 4     & 3     & 6     & 5     & 2\\
\end{array}
\right ),\ \left (
\begin{array}{cccccc}
1     & 2     & 3     & 4     & 5     & 6\\
1     & 5     & 6     & 3     & 4     & 2\\
\end{array}
\right ),\ \left (
\begin{array}{cccccc}
1     & 2     & 3     & 4     & 5     & 6\\
1     & 5     & 6     & 4     & 3     & 2\\
\end{array}
\right ),\
$$
$$
\left (
\begin{array}{cccccc}
1     & 2     & 3     & 4     & 5     & 6\\
1     & 6     & 5     & 3     & 4     & 2\\
\end{array}
\right ),\ \left (
\begin{array}{cccccc}
1     & 2     & 3     & 4     & 5     & 6\\
1     & 6     & 5     & 4     & 3     & 2\\
\end{array}
\right ),\ \left (
\begin{array}{cccccc}
1     & 2     & 3     & 4     & 5     & 6\\
2     & 3     & 4     & 5     & 6     & 1\\
\end{array}
\right ),\
$$
$$
\left (
\begin{array}{cccccc}
1     & 2     & 3     & 4     & 5     & 6\\
2     & 3     & 4     & 6     & 5     & 1\\
\end{array}
\right ),\ \left (
\begin{array}{cccccc}
1     & 2     & 3     & 4     & 5     & 6\\
2     & 4     & 3     & 5     & 6     & 1\\
\end{array}
\right ),\ \left (
\begin{array}{cccccc}
1     & 2     & 3     & 4     & 5     & 6\\
2     & 4     & 3     & 6     & 5     & 1\\
\end{array}
\right ),\
$$
$$
\left (
\begin{array}{cccccc}
1     & 2     & 3     & 4     & 5     & 6\\
2     & 5     & 6     & 3     & 4     & 1\\
\end{array}
\right ),\ \left (
\begin{array}{cccccc}
1     & 2     & 3     & 4     & 5     & 6\\
2     & 5     & 6     & 4     & 3     & 1\\
\end{array}
\right ),\ \left (
\begin{array}{cccccc}
1     & 2     & 3     & 4     & 5     & 6\\
2     & 6     & 5     & 3     & 4     & 1\\
\end{array}
\right ),\
$$
$$
\left (
\begin{array}{cccccc}
1     & 2     & 3     & 4     & 5     & 6\\
2     & 6     & 5     & 4     & 3     & 1\\
\end{array}
\right ).
$$ Let
$$P'=
\left (
\begin{array}{cccccc}
1     & 2     & 3     & 4     & 5     & 6\\
1     & 4     & 5     & 3     & 6     & 2\\
\end{array}
\right ). $$ Note that, in particular, first column of $B^{-1}$
coincide to within the permutation $P'$ with first column of
$A^{-1}$ but $P'$ is not sets isomorphism from graph $G_A$ to
graph $G_B$. One can see that $P'$ won't be such matrix after the
$3$-rd iteration of the algorithm.

{\bf Iteration 1.} Let $\varepsilon_1=0.1$; then

$ (A^1)^{-1}=\left (
\begin{array}{cccccc}
8+\varepsilon_1& 1     & 1     & 1     & 1     & 0\\
1              & 7     & 1     & 0     & 0     & 1\\
1              & 1     & 7     & 0     & 0     & 1\\
1              & 0     & 0     & 7     & 1     & 1\\
1              & 0     & 0     & 1     & 7     & 1\\
0              & 1     & 1     & 1     & 1     & 8\\
\end{array}
\right )^{-1}=$
$$=
\left (
\begin{array}{cccccc}
0.132    &-0.018    &-0.018    &-0.018    &-0.018    & 0.008811\\
-0.018   & 0.15     &-0.016    &0.004433  &0.004433  &-0.018 \\
-0.018   &-0.016    & 0.15     &0.004433  &0.004433  &-0.018\\
-0.018   &0.004433  &0.004433  & 0.15     &-0.016    &-0.018\\
-0.018   &0.004433  &0.004433  &-0.016    & 0.15     &-0.018\\
0.008811 &-0.018    &-0.018    &-0.018    &-0.018    & 0.134\\
\end{array}
\right ),$$ thus if we put $k_1=1$, then we obtain

$ (B^1)^{-1}=\left (
\begin{array}{cccccc}
8+\varepsilon_1& 0     & 1     & 1     & 1     & 1\\
0              & 8     & 1     & 1     & 1     & 1\\
1              & 1     & 7     & 1     & 0     & 0\\
1              & 1     & 1     & 7     & 0     & 0\\
1              & 1     & 0     & 0     & 7     & 1\\
1              & 1     & 0     & 0     & 1     & 7\\
\end{array}
\right )^{-1}=$

$$=
\left (
\begin{array}{cccccc}
0.132    &0.008811  &-0.018    &-0.018    &-0.018    &-0.018\\
0.008811 &0.134      &-0.018    &-0.018    &-0.018    &-0.018\\
-0.018   &-0.018    & 0.15     &-0.016    &0.004433  &0.004433\\
-0.018   &-0.018    &-0.016    & 0.15     &0.004433  &0.004433\\
-0.018   &-0.018    &0.004433  &0.004433  & 0.15     &-0.016\\
-0.018   &-0.018    &0.004433  &0.004433  &-0.016    & 0.15\\
\end{array}
\right ).$$ Correspondence $j\leftrightarrow k_j$ has the
following form after iteration 1:
$$
\left (
\begin{array}{cccccc}
1     & 2     & 3     & 4     & 5     & 6\\
1     & .     & .     & .     & .     & .\\
\end{array}
\right ),
$$
where points denote correspondences which are not settled yet.

{\bf Iteration 2.} Let us put $\varepsilon_2=0.2$

$ (A^2)^{-1}=\left (
\begin{array}{cccccc}
8+\varepsilon_1& 1              & 1     & 1     & 1     & 0\\
1              & 7+\varepsilon_2& 1     & 0     & 0     & 1\\
1              & 1              & 7     & 0     & 0     & 1\\
1              & 0              & 0     & 7     & 1     & 1\\
1              & 0              & 0     & 1     & 7     & 1\\
0              & 1              & 1     & 1     & 1     & 8\\
\end{array}
\right )^{-1}=$
$$= \left (
\begin{array}{cccccc}
0.134    &-0.017    &-0.018    &-0.018    &-0.018    & 0.00875\\
-0.017   & 0.146    &-0.016    &0.004303  &0.004303  &-0.017 \\
-0.018   &-0.016    & 0.15     &0.004447  &0.004447  &-0.018\\
-0.018   &0.004303  &0.004447  & 0.15     &-0.016    &-0.018\\
-0.018   &0.004303  &0.004447  &-0.016    & 0.15     &-0.018\\
0.00875  &-0.017    &-0.018    &-0.018    &-0.018    & 0.134\\
\end{array}
\right ),$$

$ (B^2)^{-1}=\left (
\begin{array}{cccccc}
8+\varepsilon_1& 0     & 1              & 1     & 1     & 1\\
0              & 8     & 1              & 1     & 1     & 1\\
1              & 1     & 7+\varepsilon_2& 1     & 0     & 0\\
1              & 1     & 1              & 7     & 0     & 0\\
1              & 1     & 0              & 0     & 7     & 1\\
1              & 1     & 0              & 0     & 1     & 7\\
\end{array}
\right )^{-1}=$
$$=
\left (
\begin{array}{cccccc}
0.132    &0.00875   &-0.017    &-0.018    &-0.018    &-0.018\\
0.00875  &0.134     &-0.017    &-0.018    &-0.018    &-0.018\\
-0.017   &-0.017    & 0.146    &-0.016    &0.004303  &0.004303\\
-0.018   &-0.018    &-0.016    & 0.15     &0.004447  &0.004447\\
-0.018   &-0.018    &0.004303  &0.004447  & 0.15     &-0.016\\
-0.018   &-0.018    &0.004303  &0.004447  &-0.016    & 0.15\\
\end{array}
\right ).$$ Correspondence $j\leftrightarrow k_j$ has the
following form after iteration 2:
$$
\left (
\begin{array}{cccccc}
1     & 2     & 3     & 4     & 5     & 6\\
1     & 3     & .     & .     & .     & .\\
\end{array}
\right ).
$$

{\bf Iteration 3.} $\varepsilon_3=0.3$

$ (A^3)^{-1}=\left (
\begin{array}{cccccc}
8+\varepsilon_1& 1              & 1              & 1     & 1     & 0\\
1              & 7+\varepsilon_2& 1              & 0     & 0     & 1\\
1              & 1              & 7+\varepsilon_3& 0     & 0     & 1\\
1              & 0              & 0              & 7     & 1     & 1\\
1              & 0              & 0              & 1     & 7     & 1\\
0              & 1              & 1              & 1     & 1     & 8\\
\end{array}
\right )^{-1}=$
$$=
\left (
\begin{array}{cccccc}
0.132    &-0.017    &-0.017    &-0.018    &-0.018    & 0.008659\\
-0.017   & 0.146    &-0.015    &0.004324  &0.004324  &-0.017\\
-0.017   &-0.015    & 0.144    &0.004255  &0.004255  &-0.017\\
-0.018   &0.004324  &0.004255  & 0.15     &-0.016    &-0.018\\
-0.018   &0.004324  &0.004255  &-0.016    & 0.15     &-0.018\\
0.008659 &-0.017    &-0.017    &-0.018    &-0.018    & 0.134\\
\end{array}
\right ),$$

\bigskip

$ (B^3)^{-1}=\left (
\begin{array}{cccccc}
8+\varepsilon_1& 0     & 1              & 1              & 1     & 1\\
0              & 8     & 1              & 1              & 1     & 1\\
1              & 1     & 7+\varepsilon_2& 1              & 0     & 0\\
1              & 1     & 1              & 7+\varepsilon_3& 0     & 0\\
1              & 1     & 0              & 0              & 7     & 1\\
1              & 1     & 0              & 0              & 1     & 7\\
\end{array}
\right )^{-1}=$
$$=
\left (
\begin{array}{cccccc}
0.132    &0.008659  &-0.017    &-0.017    &-0.018    &-0.018\\
0.008659 &0.134     &-0.017    &-0.017    &-0.018    &-0.018\\
-0.017   &-0.017    & 0.146    &-0.015    &0.004324  &0.004324\\
-0.017   &-0.017    &-0.015    & 0.144    &0.004255  &0.004255\\
-0.018   &-0.018    &0.004324  &0.004255  & 0.15     &-0.016\\
-0.018   &-0.018    &0.004324  &0.004255  &-0.016    & 0.15\\
\end{array}
\right ).$$ Correspondence $j\leftrightarrow k_j$ has the
following form after iteration 3:
$$
\left (
\begin{array}{cccccc}
1     & 2     & 3     & 4     & 5     & 6\\
1     & 3     & 4     & .     & .     & .\\
\end{array}
\right ).
$$

{\bf Iteration 4.} $\varepsilon_4=0.4$

$ (A^4)^{-1}=\left (
\begin{array}{cccccc}
8+\varepsilon_1     & 1                    & 1                   & 1                   & 1     & 0\\
1                   & 7+\varepsilon_2      & 1                   & 0                   & 0     & 1\\
1                   & 1                    & 7+\varepsilon_3     & 0                   & 0     & 1\\
1                   & 0                    & 0                   & 7+\varepsilon_4     & 1     & 1\\
1                   & 0                    & 0                   & 1                   & 7     & 1\\
0                   & 1                    & 1                   & 1                   & 1     & 8\\
\end{array}
\right )^{-1}=$
$$=
\left (
\begin{array}{cccccc}
0.132    &-0.017    &-0.017    &-0.017    &-0.018    & 0.008929\\
-0.017   & 0.146    &-0.015    &0.004079  &0.004351  &-0.017\\
-0.017   &-0.015    & 0.144    &0.004014  &0.004282  &-0.017\\
-0.017   &0.004079  &0.004014  & 0.142    &-0.015    &-0.017\\
-0.018   &0.004351  &0.004282  &-0.015    & 0.15     &-0.018\\
0.008541 &-0.017    &-0.017    &-0.017    &-0.018    & 0.134\\
\end{array}
\right ),$$

$ (B^4)^{-1}=\left (
\begin{array}{cccccc}
8+\varepsilon_1     & 0     & 1     & 1     & 1     & 1\\
0                   & 8     & 1     & 1     & 1     & 1\\
1                   & 1     & 7+\varepsilon_2     & 1     & 0     & 0\\
1                   & 1     & 1     & 7+\varepsilon_3     & 0     & 0\\
1                   & 1     & 0     & 0     & 7+\varepsilon_4     & 1\\
1                   & 1     & 0     & 0     & 1     & 7\\
\end{array}
\right )^{-1}=$
$$=
\left (
\begin{array}{cccccc}
0.132    &0.008541  &-0.017    &-0.017    &-0.017    &-0.018\\
0.008541 &0.134     &-0.017    &-0.017    &-0.017    &-0.018\\
-0.017   &-0.017    & 0.146    &-0.015    &0.004079  &0.004351\\
-0.017   &-0.017    &-0.015    & 0.144    &0.004014  &0.004282\\
-0.017   &-0.017    &0.004079  &0.004014  & 0.142    &-0.015\\
-0.018   &-0.018    &0.004351  &0.004282  &-0.015    & 0.15\\
\end{array}
\right ).$$ Correspondence $j\leftrightarrow k_j$ has the
following form after iteration 4:
$$
\left (
\begin{array}{cccccc}
1     & 2     & 3     & 4     & 5     & 6\\
1     & 3     & 4     & 5     & .     & .\\
\end{array}
\right ).
$$

The sets $R(A^4)$ and $R(B^4)$ are already simple before 5-th
iteration. Thus we may not perturb matrices $A^4$ and $B^4$ no
more and we can set the correspondence for the rest of the graphs
vertices ($j=\overline{5,6}$).

As a result, we obtain the isomorphism from $G_A$ to $G_B$. It is
$$
\left (
\begin{array}{cccccc}
1     & 2     & 3     & 4     & 5     & 6\\
1     & 3     & 4     & 5     & 6     & 2\\
\end{array}
\right ).
$$

\bigskip

\section{\bf COMPUTATIONAL EFFICIENCY OF THE ALGORITHM}

\subsection{Localization of separated solutions of linear equations
systems} Let us denote solution of system $A^jx=e_k$ as
$x_k^{(j)}$. By definition, put
$$R_i(A^j)=\{x_{i_1}^{(j)},\ldots, x_{i_k}^{(j)} | \forall p,q:
1\le p,q\le k\ \exists P:x_{i_p}^{(j)}=Px_{i_q}^{(j)}\}.$$ We have
$R(A)=\cup R_i(A)$. Let $m(A^j)$ be the amount of sets $R_i(A^j)$
for matrix $A^j$. $R(A)$ is simple iff $|R_i(A)|=1$ for every $i$
and hence $m(A)=n$. Let us consider these sets at the process of
the perturbations of matrix $A$.

{\bf Definition 4.} Let us say that set $R_i(A^{j-1})$ is {\it
splitted} by separation of some vector $x_j^{(j-1)}\in
R_i(A^{j-1})$ from it at $j$-th iteration of the algorithm if
$x_j^{(j)}\in R_k(A^j)$ for some $k$ and $|R_k(A^j)|=1$. Let us
call such vector as a {\it separated} vector.\qed

If we get $m(A^j)\ge m(A^{j-1})$ and, finally, $m(A^{n_0})=
m(B^{n_0})=n$ for some $n_0$ ($n_0\le n$) while perturbing the
graphs matrices at iterations of the algorithm, then sets $R(A)$
and $R(B)$ will be simple sets $R(A^{n_0})$ and $R(B^{n_0})$ at
some $n_0$-th iteration. 

By definition, put
$$\Delta_{kl}^{(j)}= |x_{k_1k_1}^{(j)}-x_{l_1l_1}^{(j)} |,\
x_{k_1}^{(j)}\in R_k(A^j),\ x_{l_1}^{(j)}\in R_l(A^j).\eqno (5)
$$
Let us understand $\Delta_{kl}^{(j)}$ as a distance between
different sets $R_k(A^j)$ and $R_l(A^j)$. If $x_{k_1k_1}^{(j)}\neq
x_{l_1l_1}^{(j)}$, then $x_{k_1}^{(j)}\in R_k(A)$ and
$x_{l_1}^{(j)}\in R_l(A)$, where $k\neq l$.

By definition, put
$$\Delta_{kl}^{(j)}= |x_{k_1k_1}^{(j)}-x_{l_1l_1}^{(j)} |,\
x_{k_1}^{(j)}\in R_k(A^j),\ x_{l_1}^{(j)}\in R_l(A^j).\eqno (5)
$$
$\Delta_{kl}^{(j)}$ may be understand as a distance between
different sets $R_k(A^j)$ and $R_l(A^j)$. If $x_{k_1k_1}^{(j)}\neq
x_{l_1l_1}^{(j)}$, then $x_{k_1}^{(j)}\in R_k(A)$ and
$x_{l_1}^{(j)}\in R_l(A)$, where $k\neq l$.

At the presented above example, the splitting of sets  may be
illustrated as follows. Before 1-st iteration of the algorithm:
$$R_1(A^0)=\{x_1^{(0)},x_6^{(0)}\},\ R_2(A^0)=\{x_2^{(0)},x_3^{(0)},x_4^{(0)},x_5^{(0)}\},$$
$$\Delta_{12}^{(0)}=0.016,$$ $$m(A^0)=2.$$ After 1-st iteration:
$$R_1(A^1)=\{x_1^{(1)}\},\ R_2(A^1)=\{x_2^{(1)},x_3^{(1)},x_4^{(1)},x_5^{(1)}\},\
R_3(A^1)=\{x_6^{(1)}\},$$ $$\Delta_{12}^{(1)}=0.018,\
\Delta_{13}^{(1)}=0.002,\ \Delta_{23}^{(1)}=0.016,$$ $$m(A^1)=3.$$
After 2-nd iteration:
$$R_1(A^2)=\{x_1^{(2)}\},\ R_2(A^2)=\{x_2^{(2)}\},\ R_3(A^2)=\{x_3^{(2)}\},$$ $$\
R_4(A^2)=\{x_4^{(2)},x_5^{(2)}\},\ R_5(A^2)=\{x_6^{(2)}\},$$
$$\Delta_{12}^{(2)}=0.014,\ \Delta_{13}^{(2)}=0.016,\
\Delta_{14}^{(2)}=0.018,\ \Delta_{15}^{(2)}=0.002,\
\Delta_{23}^{(2)}=0.004,$$ $$\Delta_{24}^{(2)}=0.004
\Delta_{25}^{(2)}=0.012,\ \Delta_{34}^{(2)}=0,\
\Delta_{35}^{(2)}=0.016,\ \Delta_{45}^{(2)}=0.016,$$
$$m(A^2)=5.$$ After 3-rd iteration:
$$R_1(A^3)=\{x_1^{(3)}\},\ R_2(A^3)=\{x_2^{(3)}\},\ R_3(A^3)=\{x_3^{(3)}\},$$ $$\
R_4(A^3)=\{x_4^{(3)},x_5^{(3)}\},\ R_5(A^3)=\{x_6^{(3)}\},$$
$$\Delta_{12}^{(3)}=0.014,\ \Delta_{13}^{(3)}=0.012,\
\Delta_{14}^{(3)}=0.018,\ \Delta_{15}^{(3)}=0.002,\
\Delta_{23}^{(3)}=0.002,$$
$$\Delta_{24}^{(3)}=0.004,\ \Delta_{25}^{(3)}=0.012,\ \Delta_{34}^{(3)}=0.006,\
\Delta_{35}^{(3)}=0.01,\ \Delta_{45}^{(3)}=0.016,$$
$$m(A^3)=5.$$ After 4-th iteration:
$$R_1(A^4)=\{x_1^{(4)}\},\ R_2(A^4)=\{x_2^{(4)}\},\ R_3(A^4)=\{x_3^{(4)}\},$$ $$\
R_4(A^4)=\{x_4^{(4)}\},\ R_5(A^4)=\{x_5^{(4)}\},\
R_6(A^4)=\{x_6^{(4)}\}.$$
$$\Delta_{12}^{(4)}=0.014,\ \Delta_{13}^{(4)}=0.012,\ \Delta_{14}^{(4)}=0.001,\
\Delta_{15}^{(4)}=0.018, \ \Delta_{16}^{(4)}=0.002,$$
$$\Delta_{23}^{(4)}=0.002,\ \Delta_{24}^{(4)}=0.004,\ \Delta_{25}^{(4)}=0.004,\ \Delta_{26}^{(4)}=0.012,\
\Delta_{34}^{(4)}=0.008,$$ $$\Delta_{35}^{(4)}=0.006,\
\Delta_{36}^{(4)}=0.01,\ \Delta_{45}^{(4)}=0.008,\
\Delta_{46}^{(4)}=0.008,\ \Delta_{56}^{(4)}=0.016,$$
$$m(A^4)=6.$$ Finally, we obtain that every initial set of
vectors $R_i(A)$ is splitted and $n_0=4$ whereas $n=6$. Note that
there are no separated vectors after 3-rd iteration but all of
column-vectors $x_j$ ($j=\overline{1,n}$) are separated from each
other after 4-th iteration.

Matrix conditional number is the following value:
$$\mu (A)=\sup\limits_{x\neq 0,\ \xi\neq 0}
\biggl\{\frac{\|Ax\|\|\xi\|}{\|A\xi\|\|x\|}\biggr\},$$ where
$\|\cdot\|$ is Euclidian norm of vectors in $\mathbb{R}^n$. Let
$\lambda_{\mbox{max}}$ be a maximal eigenvalue of matrix $A$
spectrum and let $\lambda_{\mbox{min}}$ be a minimal eigenvalue of
matrix $A$ spectrum. $$\mu (A)=\frac{\sup\limits_{x\neq
0}\|Ax\|/\|x\|} {\inf\limits_{\xi\neq 0}\|A\xi\|/\|\xi\|}=
\frac{\lambda_{\mbox{max}}}{\lambda_{\mbox{min}}}<\infty.$$
$\lambda_{\mbox{min}}\neq 0$ because $A$ is a positive definite
matrix.

Let us consider the system of equations $(A+C)y=f$ that is
obtained from $Ax=f$ by perturbation of matrix $A$ by matrix $C$.
Let
$$\frac{\|C\|}{\|A\|}\le\theta,$$
where $$\|A\|=\sup\limits_{x\neq 0}\dt\frac{\|Ax\|}{\|x\|}.$$
$\|A\|=\lambda_{\mbox{max}}$ for positive definite matrix $A$
\cite{6}. The following theorem holds \cite{6}.

{\bf Theorem (Godunov).} If $\theta\mu (A)<1$, then
$$\frac{\|y-x\|}{\|x\|}\le\frac{\theta\mu (A)}{1-\theta\mu
(A)}.$$\qed

The following holds for symmetrical matrix conditional number
\cite{6}:
$$\mu (A)\le \frac{\eta (A)}{\chi (A)},$$ where $\eta
(A)=\max\limits_{1\le i\le n}(a_{ii}+\sum\limits_{j\neq
i}|a_{ij}|),\ \chi (A)=\min\limits_{1\le i\le
n}(a_{ii}-\sum\limits_{j\neq i}|a_{ij}|).$ Since $a_{kk}=$
$=d+d_k+\varepsilon_k$ at $j$-th iteration, where $\varepsilon_k$
is value of perturbation of $k$-th matrix $A$ diagonal element and
$\varepsilon_k=0$ if $k<j$. Let $i_1$ be number of the row that
gives $\eta (A)$ and let $i_2$ be number of the row that gives
$\chi (A)$. We have
$$\eta
(A)=a_{i_1i_1}+\varepsilon_{i_1i_1}+\sum\limits_{j=1}^n|a_{i_1j}|=d+\varepsilon_{i_1i_1}+d_{i_1}+d_{i_1}=
d+\varepsilon_{i_1i_1}+d+d=3d+\varepsilon_{i_1i_1},$$ $$\chi
(A)=a_{i_2i_2}+\varepsilon_{i_2i_2}-\sum\limits_{j=1}^n|a_{i_2j}|=d+\varepsilon_{i_1i_1}+d_{i_2}-d_{i_2}=d+\varepsilon_{i_2i_2}.$$
Therefore $$\mu (A)\le \frac{\eta (A)}{\chi
(A)}=\frac{3d+\varepsilon_{i_1i_1}}{d+\varepsilon_{i_2i_2}}.$$ If
$\varepsilon_{i_1i_1}=1/n^{p_1}$,
$\varepsilon_{i_2i_2}=1/n^{p_2}$, where $p_1,p_2\in \mathbb{N}$,
then
$$\mu (A)\le
\frac{3d+\varepsilon_{i_1i_1}}{d+\varepsilon_{i_2i_2}}=\frac{3d+1/n^{p_1}}{d+1/n^{p_2}}=n^{p_2-p_1}\frac{3dn^{p_1}+1}{dn^{p_2}+1}\le
n^{p_2-p_1}\frac{4dn^{p_1}}{dn^{p_2}}=4.$$

So the stated above theorem take the following form in our case:

{\it If $\theta<1/4$, then}
$$ \|y-x\|\le\theta\frac{4}{1-4\theta}\|x\|.\eqno (6)$$
\qed

The perturbations of matrices must be small enough for guaranteed
isolation of vector $x_j^{(j-1)}$ separated from $R_i(A^{j-1})$ at
$j$-th iteration. This means that if $x_j^{(j-1)}$ is separating
from $R_i(A^{(j-1)})$, then it is necessary that, after the
perturbation, $x_j^{(j-1)}$ will be sufficiently far from other
sets in sense of introduced distance $\Delta$ (5). We can
guarantee this if there is a localization of the linear equations
systems solutions take place. Here, we means localization about
the splitted set $R_i(A^{(j-1)})$ in the sense of distance
$\Delta$.

It follows from (6) that
$$\forall i:\ \|x_i^{(j)}-x_i^{(j-1)}\|\le\frac{\theta(\varepsilon_j)4}{1-4\theta(\varepsilon_j)}\|x_i^{(j-1)}\|,$$ where
$$\frac{\|C\|}{\|A^{j-1}\|}=\frac{\|\varepsilon_jE^j\|}{\|A^{j-1}\|}
=\frac{\varepsilon_j}{\|A^{j-1}\|}.$$ So we can put
$$\theta(\varepsilon_j)=\dt\frac{\varepsilon_j}{\|A^{j-1}\|}.$$ Putting $\varepsilon_j=1/n^p$, we get
$$\|x_i^{(j-1)}-x_i^{(j)}\|<
\frac{4/n^p}{\|A^j\|-4/n^p}\|x_i^{(j-1)}\|=\frac{4\cdot
1/n^p}{\|A^j\|-4\cdot 1/n^p}\|x_i^{(j-1)}\|=$$$$=
\frac{4}{\|A^j\|n^p-4}\|x_i^{(j-1)}\|
=\frac{4}{\|A^j\|n^p+o(n)}\|x_i^{(j-1)}\|\le
\frac{4}{\|A^j\|n^p}\|x_i^{(j-1)}\|.\eqno (7)$$

We have
$$\mu(A^{j-1})=\lambda_{\mbox{max}}/\lambda_{\mbox{min}}\le
4,$$ where $\lambda_{\mbox{max}}$ is a maximal eigenvalue of
matrix $A^{j-1}$ spectrum and $\lambda_{\mbox{min}}$ is a minimal
eigenvalue of matrix $A^{j-1}$ spectrum. Since Gershgorin's
Theorem \cite{7}, we have
$$\lambda_{\mbox{max}}\ge d.$$ Taking into account
that $\lambda_{\mbox{min}}<\|A^{j-1}\|$, we get
$$\|A^{j-1}\|>\lambda_{\mbox{min}}>\lambda_{\mbox{max}}/4>d/4.$$
So we have
$$\|x_i^{(j-1)}-x_i^{(j)}\|<\frac{16}{dn^p}\|x_i^{(j-1)}\|.$$
Since $\lambda_{\mbox{min}}>1$, $\|e_j\|=1$ and
$\|x_i^{(j-1)}\|=A^{-1}e_j$, then $\|x_i^{(j-1)}\|<1$. So we have
$$\|x_i^{(j-1)}\|<\frac{\sqrt{n}}{d}.$$ Therefore it follows from (7) that  $$\forall i:\ \|x_i^{(j-1)}-x_i^{(j)}\|
<\frac{4}{d^2n^{p-1}}.$$ It is easy to prove that
$$|x_{ii}^{(j-1)}-x_{ii}^{(j)}|\le\|x_i^{(j-1)}-x_i^{(j)}\|.$$ So we have
$$\forall i:\ |x_{ii}^{(j-1)}-x_{ii}^{(j)}|<\frac{4}{d^2n^{p-1}}\eqno (8)$$
(8) gives an estimate of localization of vectors from sets
$R_i(A^j)$ at the process of the algorithm iteration.

\subsection{Splitting of the sets of linear equations systems
solutions}

Implemented perturbations must be sufficient for splitting of the
sets $R_k$ at the process of operation of the algorithm. The
algorithm is computationally effective only if,  for every
$j,k,l$, distances $\Delta_{kl}^{(j)}$ will be a non-zero machine
number. Let us show that we can obtain adequate accuracy using
machine numbers which mantissas length not exceeding $n$.

Let $x_i^{(j-1)},\ x_j^{(j-1)}\in  R_i(A^{j-1})$ for some $i$,
i.e., there are at least two elements that prevents us to set the
univocal correspondence between the vertices of graphs at $j$-th
iteration. If $x_{ii}^{(j)}\neq x_{jj}^{(j)}$ after $j$-th
iteration, then there is no permutation matrix $P$ such that
$x_j^{(j)}=Px_i^{(j)}$. Let us estimate the lower bound of the
distance between such $x_i^{(j)}$ and $x_j^{(j)}$ after $j$-th
iteration of the algorithm.

{\bf Theorem 1.} If there exists permutation matrix $P$ such that
$x_j^{(j-1)}=Px_i^{(j-1)}$ and $0<\varepsilon_j<1$, then
$$\Delta_{ij}^{(j)}=|x_{jj}^{(j)}-x_{ii}^{(j)}|>\frac{1}{3^nd^2(3d/\varepsilon_j+1)}.$$

{\bf Proof.} Put $A=A^{j-1},\ A'=A^j,\ x_i=x_i^{(j-1)},\
x_i'=x_i^{(j)},\ x_j=x_j^{(j-1)},\ x_j'=x_j^{(j)}$. If
$x_j'=Px_i'$, then $x_{jj}'=x_{ii}'$.
$$x_{jj}'=\dt\frac{A_{jj}'}{\det A'},\ x_{ii}'=\dt\frac{A_{ii}'}{\det A'}.$$ We have
$$\det A'=\det A+\varepsilon_jA_{jj},\ A_{jj}'=A_{jj}$$ because the diagonal element $a_{jj}$ is the
only element of $A$ that changes at the iteration. Since
$$A_{ii}'=A_{ii}+\varepsilon_jA_{ii,jj},$$ we get
$$|x_{jj}'-x_{ii}'|=\biggl |\dt\frac{A_{jj}'}{\det
A'}-\dt\frac{A_{ii}'}{\det A'} \biggr|=\biggl
|\dt\frac{A_{jj}'}{\det
A'}-\dt\frac{A_{ii}+\varepsilon_jA_{ii,jj}}{\det A'}
\biggr|=\dt\frac{\varepsilon_jA_{ii,jj}}{\det
A+\varepsilon_jA_{jj}}.\eqno(9)$$

$A_{ii,jj}$ is a determinant of the submatrix of $A$ that is
obtained by deleting of $i$-th row and $i$-th column and $j$-th
row and $j$-th column. $A$ is a symmetrical positive definite
matrix with a diagonal predominance hence the following holds for
$A$ (Hadamard):
$$H_k\equiv |a_{kk}|-\sum\limits_{l\neq k}|a_{kl}|=d_{k}>0,\ i=\overline{1,n},$$
where $$ d_{k}'=\left\{
\begin{array}{ll}
d+\varepsilon_k,&\mbox{ if } 1\le k\le j,\\
d,&\mbox{ if } j< k\le n.\\
\end{array}
\right.
$$
Consequently (see \cite{8}), $$A_{ii,jj}\ge H_1\ldots
\widehat{H_i}\ldots\widehat{H_j}\ldots H_n=d_1'\ldots
\widehat{d_i'}\ldots\widehat{d_j'}\ldots d_n'>d^{n-2},\eqno(10)$$
where $\widehat{H_i}$ means that $H_i$ is absent at the product.

On the other hand, there is a following estimate for the maximal
eigenvalue $\lambda_{\mbox{max}}$ of $A$ (Gershgorin's Theorem,
see \cite{7}):
$$\lambda_{\mbox{max}}\le
3d+\max\limits_{1\le k\le j}\varepsilon_k<3d+1,$$ if $\forall k\
\varepsilon_k<1$. Therefore
$$\det A=\prod_{k=1}^n\lambda_k<\lambda_{\mbox{max}}^n<(3d+1)^n,\eqno(11)$$

$$A_{jj}=\prod_{k=1}^{n-1}\lambda_k'<\lambda_{\mbox{max}}^{n-1}<(3d+1)^{n-1},\eqno(12)$$
Let $\lambda_{\mbox{max}}'$ be the maximal eigenvalue of the
matrix that is remaining from the matrix $A$ by deleting its
$i$-th and $j$-th rows and columns. $\lambda_{\mbox{max}}'\le
\lambda_{\mbox{max}}$ \cite{8}, consequently, combining (10),
(11), (12), from (9) we get:
$$|x_{jj}'-x_{ii}'|=\dt\frac{\varepsilon_jA_{ii,jj}}{\det A+\varepsilon_jA_{jj}}>
\dt\frac{\varepsilon_jd^{n-2}}{(3d+1)^n+\varepsilon_j(3d+1)^{n-1}}>\dt\frac{\varepsilon_jd^{n-2}}{(3d)^{n+1}+\varepsilon_j(3d)^n},$$
i.e.,
$$\Delta_{ij}^{(j)}=|x_{jj}^{(j)}-x_{ii}^{(j)}|>\dt\frac{1}{3^nd^2(3d/\varepsilon_j+1)}.\eqno(13)$$
\qed

\bigskip

Let we shall be solve the systems of linear equations $A^jx=e_k$
and $B^jy=e_k$ ($j,k=\overline{1,n}$) by Gauss-Seidel iterative
method. Actually, it is doesn't really matter what method of
solution is chosen: we just need a proper accuracy of solutions
which we can obtain in a proper time. But Gauss-Seidel method also
gives us an efficient way to estimate the needed machine number
mantissas length.

The following theorem may be found in \cite{9}.

{\bf Theorem 2}. Let $\sum\limits_{j\neq i}|a_{ij}|\le \gamma
|a_{ii}|, \gamma<1,$ for every $i=\overline{1,n}$. Then
$$\|x_i-\tilde{x}_i^{s}\|\le\gamma\|x_i-\tilde{x}_i^{s-1}\|,$$
where $x_i$ is exact solution of the system of linear equations
$A^jx=e_i$, $\tilde{x}_i^{s}$ is an approximate solution of the
system of equations at $s$-th iteration of Gauss-Seidel
method.\qed

If $A=A^j$ for some $j$, then $\gamma\le 1/2$. Therefore at $s$-th
iteration of Gauss-Seidel method we have:
$$|x_{ii}-\widetilde{x}_{ii}^s|\le\|x_i-\widetilde{x}_i^s\|<\frac{1}{2^s}\delta^0,\
|x_{jj}-\widetilde{x}_{jj}^s|\le\|x_j-\widetilde{x}_{j}^s\|<\frac{1}{2^s}\delta^0,$$
where $\delta_0$ is a mistake of initial approximation. If
$\varepsilon_j=1/n^p$ at $j$-th iteration, then, taking into
account (12), we get
$$|x_{jj}^{(j)}-x_{ii}^{(j)}|>\frac{1}{3^nd^2(3dn^p+1)}.$$
Hence if
$$\frac{1}{2^{s-1}}\delta^0\le\frac{1}{3^nd^2(3dn^p+1)},\eqno (14)$$
then difference $|x_{jj}^{(jj)}-x_{ii}^{(ii)}|$ may be fixed by
machine numbers which mantissas length is restricted by $n$. Let
us estimate the number of Gauss-Seidel method iterations $s$ that
must be implemented to fix this value. Using (14), we get
$$3^nd^2(3dn^p+1)\delta^0\le 2^{s-1}.\eqno(15)$$ Taking the
logarithm of (15), we obtain
$$n\log_23+\log_2d^2(3dn^p+1)+\log_2\delta^0\le s-1,$$ i.e.,
the number of needed method's iterations may be estimated as
$$s\ge n\log_23+p\log_2n+3\log_2d+\log_2\delta^0+3.\eqno(16)$$

Let $p>0$ be an integer that defines the value of matrix diagonal
element perturbation: $\varepsilon_j=1/n^p$. Let $x_j^{(j-1)}$ be
a vector that isolates from set $R_i(A^{j-1})$ at $j$-th
iteration.

Let us choose $p>0$ such that
$$\frac{4}{d^2n^p}<\frac{\Delta_j}{n},$$ i.e.,
$$p>\log_n\frac{4}{d^2\Delta_j}+1,\eqno (17)$$ where $\Delta_j=\min\limits_{k\neq
l}\Delta_{kl}^{(j)}$. It follows from (7) that
$$\forall i:\ |x_{ii}^{(j-1)}-x_{ii}^{(j)}|<\frac{4}{d^2n^p}<\dt\frac{\Delta_j}{n}.$$ And, using of theorem 1, we have
$$\forall i:\ \dt\frac{1}{3^nd^2(3d/\varepsilon_j+1)}<|x_{jj}^{(j)}-x_{ii}^{(j)}|<\frac{4}{d^2n^p}$$
after $j$-th iteration. This means that, at $j$-th iteration,  we
can choose such perturbation that will give us a separation that
is fixable by machine numbers with restricted mantissas length and
localization of vectors $x_i$ at the same time. This mean that we
can set the interval in which the separating of vectors occurs. We
can perform this process in the way such that if we isolating some
vector $x_j$ from $R_i(A^{j-1})$, then this vector won't be placed
in another such set of vectors at the next iteration unless power
of this set is equal to 1.

{\bf Lemma 4.} $$\Delta_j>\dt\frac{1}{3^nd^2(3dn^n+1)}$$ and
$$s_j\le (n+3)\log_2n$$ at every $j$-th iteration of the algorithm.

{\bf Proof.} Let us consider the following case of the graph
isomorphism problem. Let graphs $G_A$ and $G_B$ be complete
graphs, i.e., every vertex of each graph is adjacent to every
other vertex of the graph.  It is clear that checking of the
isomorphism of complete graphs is easy taken by itself since every
complete graph is isomorphic to another complete graph on the same
vertices set. The purpose of introducing the following example is
to show the implementation of the algorithm in hardest case
because $m(A)=1$, $|R_1(A)|=n$ and $d=n$ here, i.e $d$ is maximal
for every non-weighted graph, so, taking into account theorem 1,
$\Delta_{j}$ take on a minimal value in this case.

\noindent 1. Put $\varepsilon_1=1/n$ at 1-st iteration.
$x_1^{(1)}$ is separated at this iteration and, for every $j\neq
1$, $$1/(3^nd^2(3dn+1))\le|x_{11}-x_{jj}|\le 4/(d^2n).$$ Needed
number of Gauss-Seidel method iterations is $s$, where
$$s=n\log_23+\log_2n+3\log_2d+\log_2\delta^0+3.$$

\noindent 2. Put $\varepsilon_2=1/n^2$ at 2-nd iteration. 2-nd
vector is separated at this iteration and, for every $j\neq 2$,
$$1/(3^nd^2(3dn^2+1))\le|x_{22}-x_{jj}|\le\ 4/(d^2n^2).$$ Needed
number of Gauss-Seidel method iterations is $s$, where
$$s=n\log_23+2\log_2n+3\log_2d+\log_2\delta^0+3.$$

\noindent Continuing separating solutions of linear equations
systems in this way, we obtain $\varepsilon_2=1/n^{n-1}$ at
$(n-1)$-th iteration of the algorithm and $(n-1)$-th vector is
separated and, for every $j\neq n-1$,
$$1/(3^nd^2(3dn^{n-1}+1))\le|x_{n-1,n-1}-x_{jj}|\le
4/(d^2n^{n-1}).$$ Needed number of Gauss-Seidel method iterations
is $s$, where
$$s=n\log_23+(n-1)\log_2n+3\log_2d+\log_2\delta^0+3.$$

Since exact solution of the systems belongs to segment
$[0,1]^n\in\mathbb{R}^n$, we may put $\delta^0=O(\sqrt n)$ and,
consequently, $s\le (n+3)\log_2n$. \qed

As a result, if we set machine numbers mantissas length is equal
to $n$ and perform $s=(n+3)\log_2n$ Gauss-Seidel method iterations
for solution of every system of linear equations at iterations of
the algorithm, then we obtain that $m(A^{n-1})=n$ while the
algorithm is computationally effective to fix occurring split of
sets of linear equations systems solutions. If we use the {\it
extended} numerical machine numbers type that realized in
languages such as Object Pascal and C++, then we can operate with
machine numbers at the range from $3,6\times 10^{-4951}$ to
$1,1\times 10^{4932}$. Thus the algorithm may be efficient for
reconstructible graphs for which the number of vertices is lesser
than 4951 and it is possible to operate with standard extended
machine numbers without any additional procedures for correct
operating with them.

\subsection{Overall complexity of the algorithm}

At every $j$-th iteration of the algorithm, we checks similarity
of $G_A^{1,\ldots ,j}$ and $G_B^{k_1,\ldots ,k_{j-1},l}$ at step
3, where $l$ can take on a value from 1 to $n$. Checking of
similarity implements using inverse matrices to matrices $A^j$ and
$B^l$ ($l=\overline{1,n}$), where
$$B^l=B+\sum\limits_{i=1}^{j-1}\varepsilon_iE^{k_i}+\varepsilon_iE^l.$$
We obtain inverse matrices solving the systems of linear
equations. Thus we must solve $n+n\times n$ linear equations
systems at most. Let $N_S$ be the number of elementary machine
operations that we must implement to solve linear equations system
with needed accuracy. Then complexity of finding all needed
inverse matrices is equal to $n(n+1) N_S$.

Checking of coincidence to within the permutations of inverse
matrices columns can be implemented at $O(n\log n)$ elementary
machine operations. To check the similarity of two graphs
$G_A^{1,\ldots ,j}$ and $G_B^{k_1,\ldots ,k_{j-1},l}$ we must
implement checking of $n^2$ pairs of inverse matrices columns at
most. Therefore complexity of checking for similarity of two
graphs is equal to $O(n^3\log n)$. The number of such graph pairs
may be equal to $n$ at most.

Hence overall complexity of one iteration of the algorithm is
equal to $$O(n(n+1)N_S)+n\cdot O(n^3\log n)=O(n^2N_S+n^2\log
n).\eqno (18)$$

Let us estimate the value of $N_S$ for Gauss-Seidel method of
solution of linear algebraic equations.

At $j$-th iteration of the algorithm, we calculate $p$ such that
$4/(d^2n^p)<\Delta_j/n$, where $\Delta_j=\min\limits_{k\neq
k}\Delta_{kl}^{(j)}$, i.e.,
$$p=\log_n\frac{4}{d^2\Delta_j}+1.$$ Then we implement $s$
iterations of Gauss-Seidel method, where, from lemma 4, we have
$$s=(n+3)\log_2n.$$ Let $N_{S^j}$
be a complexity linear equations systems solution at $j$-th
iteration. Then
$$N_{S}^j=O(n^2)\cdot(n+3)\log_2n$$
since one iteration of Gauss-Seidel method has complexity that is
equal to $O(n^2)$. Thus, using (18), summary complexity of all
iterations of the algorithm may be estimated as
$$O (n^2)\cdot O(n^2)\cdot(n+3)\log_2n=O(n^6+n^5\log_2n).\eqno(19)$$ 

Complexity of the algorithm may be reduced if the sets
$R(A^{n_0})$ and $R(B^{n_0})$ become simple at $n_0$-th iteration,
where $n_0$ is much lesser than $n$. In particular, applying the
algorithm for isomorphism checking of regular graphs which are
lattices on torus, the sets $R(A^{n_0})$ and $R(B^{n_0})$ become
simple when $n_0=\sqrt{n}$. Vertex degree is equal to 4 at this
case. Let us remark that regular graph isomorphism testing is one
of the hardest case of the graph isomorphism problem for most
efficient algorithms which are designed for solving of the
problem. As an example, complexity of NAUTY algorithm \cite{10}
becomes exponential when vertex degree of regular graphs is equal
to 4 \cite{11}.

{\bf Remark 2.} Since (11) and (12) are obtained using
majorization of graph matrix $A$ eigenvalues by
$\lambda_{\mbox{max}}$, then (13) holds regardless of multiplicity
of form (1) matrix eigenvalues and (13) holds regardless of
multiplicity of graph adjacency matrix eigenvalues too. Thus (20)
holds regardless of multiplicity of of reconstructible graphs
adjacency matrix eigenvalues.

{\bf Remark 3.} Since $d\le n$, we have
$$|x_{jj}'-x_{ii}'|>\frac{1}{3^nd^2(3d/\varepsilon_j+1)}\ge
\frac{1}{3^nn^2(3n/\varepsilon_j+1)}.$$  Thus (20) holds
regardless of maximal degree of reconstructible graphs vertices.

\section*{\bf CONCLUSIONS}

An algorithm for solving of the graph isomorphism problem is
presented. The algorithm solve the problem for the class of
reconstructible graphs which is defines at the paper. The
algorithm checks weather graphs belongs to the class of
reconstructible graphs or not at the process of algorithm
implementing.

Proved the theorem that gives an estimate of computational
complexity of the algorithm. It is shown that algorithm is
polynomial in sense of using elementary machine operations and in
the sense of used memory too since used machine numbers mantissas
length is restricted by polynomial of number of graph vertices. It
is shown that this is holds regardless of maximal degree of the
graphs, graphs genus, graph eigenvalue multiplicity etc., i.e.,
using the algorithm, solution of the graph isomorphism problem has
no specific that may be determinated by any graph characteristics
that is usually considered.

\end{document}